\documentclass[12pt,reqno,a4paper,twoside]{article}

\usepackage{amsmath,amsthm,amstext,amscd,amssymb,euscript,mathrsfs}
\usepackage{epsfig}
\usepackage{color}
\usepackage{comment}

\usepackage{dsfont}

\newcommand{\dx}{\textup{d}x}

\newcommand{\Z}{\mathbb Z}
\newcommand{\R}{\mathbb R}
\newcommand{\N}{\mathbb N}

\newcommand{\E}{\mathbb E}

\renewcommand{\Pr}{\mathbb P}

\renewcommand{\phi}{\varphi}

\newcommand{\La}{\ensuremath{\Lambda}}

\newcommand{\si}{\ensuremath{\sigma}}

\newcommand{\pee}{\ensuremath{\mathbb{P}}}

\newcommand{\dd}{\textup{d}}

\newtheorem{theorem}{{\small T}{\scriptsize HEOREM}}[section]
\newtheorem{corollary}{{\bf{\small C}{\scriptsize OROLLARY}}}[section]
\newtheorem{proposition}{{\bf{\small P}{\scriptsize ROPOSITION}}}[section]
\newtheorem{lemma}{{\bf{\small L}{\scriptsize EMMA}}}[section]
\newtheorem{remark}{{\bf{\small R}{\scriptsize EMARK}}}[section]
\newtheorem{definition}{{\bf{\small D}{\scriptsize EFINITION}}}[section]

\renewenvironment{proof}[1]
{\noindent{{\bf{\small{ P}{\scriptsize ROOF}}}.}\hspace{0.1cm} #1} {$\;\qed$\newline}

\newcommand{\beq}{\begin{eqnarray}}
\newcommand{\eeq}{\end{eqnarray}}

\newcommand{\ba}{\begin{align*}}
\newcommand{\ea}{\end{align*}}

\newcommand{\be}{\begin{equation}}
\newcommand{\ee}{\end{equation}}

\newcommand{\bl}{\begin{lemma}}
\newcommand{\el}{\end{lemma}}

\newcommand{\br}{\begin{remark}}
\newcommand{\er}{\end{remark}}

\newcommand{\bt}{\begin{theorem}}
\newcommand{\et}{\end{theorem}}

\newcommand{\bd}{\begin{definition}}
\newcommand{\ed}{\end{definition}}

\newcommand{\bp}{\begin{proposition}}
\newcommand{\ep}{\end{proposition}}

\newcommand{\bc}{\begin{corollary}}
\newcommand{\ec}{\end{corollary}}

\newcommand{\bpr}{\begin{proof}}
\newcommand{\epr}{\end{proof}}

\newcommand{\bi}{\begin{itemize}}
\newcommand{\ei}{\end{itemize}}

\newcommand{\ben}{\begin{enumerate}}
\newcommand{\een}{\end{enumerate}}

\newcommand{\caA}{{\mathcal A}}

\newcommand{\caF}{{\mathcal F}}
\newcommand{\caG}{{\mathcal G}}

\newcommand{\caK}{{\mathcal K}}

\newcommand{\caM}{{\mathcal M}}

\newcommand{\caT}{{\mathcal T}}

\newcommand{\caZ}{{\mathcal Z}}

\DeclareMathOperator{\lawto}{\stackrel{\scriptscriptstyle\rm law}{\longrightarrow}}

\def\1{\mathds{1}} 

\begin{document}
\title{Nonconventional averages\\
along arithmetic progressions\\ and lattice spin systems}
\author{
Gioia Carinci$^{\textup{{\tiny(a)}}}$,
Jean-Ren\'e Chazottes$^{\textup{{\tiny(b)}}}$,\\
Cristian Giardin{\`a}$^{\textup{{\tiny(a)}}}$,
Frank Redig$^{\textup{{\tiny(c)}}}$.\\
{\small $^{\textup{(a)}}$
University of Modena and Reggio Emilia}\\
{\small via G. Campi 213/b, 41125 Modena, Italy}
\\
{\small $^{\textup{(b)}}$ Centre de Physique Th\'eorique, CNRS, \'Ecole Polytechnique}\\
{\small 91128 Palaiseau Cedex, France}\\
{\small $^{\textup{(c)}}$ Delft Institute of Applied Mathematics, Technische Universiteit Delft}\\
{\small Mekelweg 4, 2628 CD Delft, The Netherlands}\\
}
\maketitle

\begin{abstract}
We study the so-called nonconventional averages in the context of
lattice spin systems, or equivalently random colourings of the integers.
For i.i.d. colourings, we prove a large deviation principle
for the number of monochromatic arithmetic progressions of size two
in the box $[1,N]\cap \N$, as $N\to\infty$,  with an explicit rate function
related to the one-dimensional Ising model.\newline
For more general colourings, we prove some bounds
for the number of monochromatic arithmetic progressions
of arbitrary size, as well as for the maximal
progression inside the box $[1,N]\cap \N$.\newline
Finally, we relate nonconventional sums along arithmetic progressions
of size greater than two to statistical mechanics models in dimension
larger than one.
\end{abstract}


\section{Introduction}

Nonconventional averages along arithmetic progressions
are averages of the type
\be\label{boeloemba}
\frac1N\sum_{i=1}^N f_1 (X_{i})f_2 (X_{2 i})\cdots f_k (X_{\ell i})
\ee
where $(X_n)$ is a sequence of random variables,
and $f_1,\ldots, f_\ell$ are bounded measurable functions.

Motivation to study such averages comes from the study of arithmetic progressions
in subsets of the integers, and multiple recurrence
and multiple ergodic averages. In that context,
typically $X_i= T^i (x)$, with $T$ a weakly mixing transformation, and
$x$ is distributed according to the unique invariant measure.
See e.g. \cite{furst,mat,kra} for more background on this deep and growing field.

Only recently, starting with the work of Kifer
\cite{k}, and Kifer and Varadhan \cite{kv}, central limit behavior of
nonconventional averages was considered. These authors consider
averages along progressions more general than the arithmetic ones.
It is natural to consider the averages of
the type \eqref{boeloemba} from a probabilistic point of view and
ask questions such as whether they satisfy a large deviation principle,
whether associated extremes have classical extreme value behavior, etc.

These questions are far from obvious, since even in the simplest case
of $f_i$ being all identical, the sum
\[
S_N=\sum_{i=1}^N \prod_{j=1}^\ell f(X_{ji})
\]
is quite far from a sum of shifts of a local function.
In particular it is highly non-translation invariant.
From the point of view of statistical mechanics, large deviations
of $S_N/N$ are related to partition function and free energy
associated to the ``Hamiltonian'' $S_N$. Since $S_N$ is not translation-invariant
and (extremely) long-range, even the existence of the associated free energy is
not obvious.

In this paper, we restrict to random variables $X_i$ taking
values in a finite set.
For the sake of definiteness, we assume the joint distribution to be
a Gibbs measure with an exponentially decaying interaction to obtain
fluctuation properties of $S_N$ in a straightforward way. 
In Section \ref{basicproperties} we obtain some basic probabilistic properties
using Gaussian concentration and Poincar\'e's inequality which are available
for the Gibbs measures we consider.\newline
In Section \ref{sect-size2}, we explicitly compute the large deviation
rate function of $\frac1N\sum_{i=1}^N X_i X_{2i}$
when the $X_i$'s are  i.i.d. Bernoulli random variables.
Even if this is the absolute simplest setting, the rate function turns out to be an interesting
non-trivial object related to the one-dimensional Ising model.
Recently there has been a lot of interest in multifractal analysis
of non-conventional ergodic averages 
\cite{kps1,kps2,PS,FLM,FSW}.
Large deviation rate functions are often related 
to multifractal spectra of conventional ergodic averages.
In the present context,  this connection is not as straightforward as it is in the
context of sums of shifts of a local function.
We expect the results of this paper to be useful
in establishing such connection in the context
of non-conventional averages.

Finally, we analyze in the last section the case of arithmetic progressions
of size larger than two. This  naturally leads to statistical 
mechanics models in dimension higher than one, possibly having
phenomenon of phase transitions. Conversely the classical
Ising model in dimension $d>1$ can be related to specific
unconventional sums, which we describe below. Such a connection
deserves future investigations.

\section{The setting}

We consider $K$-colorings of the integers and denote them as $\si,\eta,$ elements of
the set of configurations
$\Omega=\{0,\ldots,K\}^{\Z}$.
We assume that on $\Omega$ there is a translation-invariant
Gibbs measure with an exponentially decaying interaction, denoted by $\pee$.
This means that, given $\alpha\in \{0,\ldots,K\}$, for the one-site conditional probability
\[
\phi_{\hat{\si}}(\alpha)=\pee(\si_0=\alpha|\si_{\Z\setminus\{0\}}= \hat{\si})
\]
we assume the variation bound
\[
\|\phi_{\hat{\si}}-\phi_{\hat{\eta}}\|_\infty
\leq e^{-n\rho}
\]
for some $\rho>0$ whenever $\hat{\si}$ and $\hat{\eta}$ agree on $[-n,n]\cap\{ {\Z\setminus\{0\}}\}$.
This class of measures is closed under single-site transformations, i.e., if we
define new spins $\si'_i = F(\si_i)$ with $F:\{0,1,\ldots,K\}\to \{0,1\ldots, K'\}$, $K'<K$,
then $\pee'$, the image measure on $\{0,1\ldots, K'\}^{\Z}$, is again a Gibbs measure
with exponentially decaying interaction, see e.g.\ \cite{red} for a proof.
In the last section,
we restrict to product measures.

For the rest of the paper we consider only $2$-colorings (i.e. $K=1$).
Given an integer $\ell$, we are interested in the random variable
\[
\sum_{i=1}^N \prod_{j=1}^\ell \si_{ji}
\]
which counts the number of arithmetic progressions of size $\ell$
with ``colour'' $1$ (starting from one) in the block $[1,Nk]$.



If we consider $K$-colorings and monochromatic arithmetic progressions, i.e.,
random variables of the type
\[
\sum_{i=1}^N \prod_{j=1}^\ell \1(\si_{ji}=\alpha)
\]
for given $\alpha\in \{0,\ldots,K\}$, then we can define the new ``colors''
$\si'_i= \1(\si_i=\alpha)$ which are zero-one
valued and, as stated before, are distributed according to $\pee'$,
a Gibbs measure with an exponentially decaying interaction.
Therefore, if we restrict to monochromatic arithmetic progressions, there is
no loss of generality if we consider $2$-colorings.

Define the averages
\[
\caA^\ell_N =\frac1N \sum_{i=1}^N \prod_{j=1}^\ell \si_{ji}.
\]
Several natural questions can be asked about them
and about some related quantities. 
We give here a non-exhaustive list. Questions 1 and 2 on this list
have been answered positively in the literature in a much
more general context (see \cite{furst} for question 1 and \cite{k, kv}
for question 2). On the contrary questions 3 and 4 have not
been considered before.
\ben
\item {\em Law of large numbers}: Does $\caA^\ell_N$ converge to $(\E (\si_0))^\ell$ as $N\to\infty$ with
$\pee$ probability one ?
\item {\em Central limit theorem}: Does there exist some $a^2>0$ such that
\[
\sqrt{N}\big(\caA^\ell_N-\left(\E (\si_0)\right)^\ell\big)\lawto {\cal N}(0,a^2),\, \textup{as}\, N\to\infty \, ?
\]
\item {\em Large deviations}: Does the rate function
\[
I(x)=\lim_{\epsilon\to 0}\lim_{N\to\infty} -\frac{1}{N} \log \pee\left(\caA^\ell_N\in[x-\epsilon,x+\epsilon]\right)
\]
exist and have nice properties ?
In view of the G\"{a}rtner-Ellis theorem \cite{dz}, the natural candidate for $I$ is
the Legendre transform of the ``free-energy''
\[
\caF(\lambda)= \lim_{N\to\infty}\frac1N \log \E\big(e^{N\lambda\caA^\ell_N}\big)
\]
provided this limit exists and is differentiable.
If, additionally, $\caF$ is analytic in a neighborhood of the origin, then
the central limit theorem follows \cite{bric}.
\item {\em Statistics of nonconventional patterns.} Let
\begin{align*}
 & \caM(N)= \\ 
& \qquad \max\{ k\in\N: \exists \; 1\leq i \leq N/k\ \text{such that}\ \si_i=1,
\si_{2i}=1,\ldots, \si_{ki}=1\}
\end{align*}
be the maximal arithmetic progression  of colour $1$ starting from zero
in the block $[1,N]$. One would expect
\[
\caM(N)\approx C \log N + X_N
\]
where $0<C<\infty$ and $X_N$ is a tight sequence of random variables with an approximate Gumbel
distribution, i.e.,
\[
e^{-c_1 e^{-x}} \leq \pee \left(X_N\leq x\right)\leq e^{-c_2 e^{-x}}.
\]
Related to this is the exponential law for the occurence
of ``rare arithmetic progressions'':
Let
\[
\caT(\ell)=\inf \{ n\in \N: \exists \; 1\leq i\leq n/\ell\ \text{such that}\ \si_i=1, \si_{2i}=1,\ldots, \si_{\ell i}=1\}
\]
be the smallest block $[1,n]$ in which a monochromatic arithmetic progression can be found
with size $k$. Then one expects that $\caT(\ell)$, appropriately
normalized, has approximately (as $\ell\to\infty$) an
exponential distribution.
Finally, another convenient quantity is
\[
\caK(N,\ell)= \sum_{i=1}^{\lfloor N/\ell\rfloor} \prod_{j=1}^\ell \si_{ji}=
\lfloor N/\ell\rfloor \caA^\ell_{\lfloor N/\ell\rfloor}
\]
which counts the number of monochromatic arithmetic progressions of size $\ell$ inside
$[1,N]$.\newline
The probability distributions of these quantities are related
by the following relations:
\[
\Pr(\caK(N,\ell) = 0) =
\Pr(\caM(N)   < \ell) =
\Pr(\caT(\ell)   > N) \;.
\]

\een

\section{Some basic probabilistic properties}\label{basicproperties}

In this section we prove some basic facts about the nonconventional averages considered
in the previous section.
\bp
\ben
\par\noindent
\item {\em Gaussian concentration bound}. Let $\ell\geq 1$ be an integer.
There exists a constant $C>0$ such that for all $n\geq 1$ and all $t>0$
\be\label{colibri}
\pee\left(\left|\caA_N^\ell-\E\left(\caA_N^\ell\right)\right|>t\right)\leq e^{-CN t^2}.
\ee
In particular, $\caA_N^\ell$ converges almost surely to $\left(\E (\si_0)\right)^\ell$
as $N$ goes to infinity.
\item {\em Logarithmic upper bound for maximal monochromatic progressions}.
There exists $\gamma>0$ such that for all $c>\gamma$
\[
\caK (N, c\log N)\to 0
\]
in probability as $N\to\infty$.
\een
\ep
\bpr
A Gibbs measure for an exponentially decaying interaction satisfies both the Gaussian concentration
bound (see e.g. \cite{chaz2}), and the Poincar\'{e} inequality \cite{chaz}.
For  a bounded measurable function $f:\Omega\to\R$ let
\[
\nabla_i f(\si) = f(\si^i)-f(\si)
\]
be the discrete derivative at $i\in\Z$, where $\si^i$ is the configuration obtained from $\si\in\Omega$
by flipping the symbol at $i$.
Next define the variation
\[
\delta_i  f=\sup_\si \nabla_i f(\si)
\]
and
\[
\|\delta f\|_2^2 =\sum_{i\in\Z}(\delta_i f)^2.
\]
Then, on the one hand,  we have the Gaussian concentration inequality: there exists some $c_1>0$ such
that
\be\label{gauss}
\pee\left(|f-\E(f)|>t\right)\leq e^{-\frac{c_{1} t^2}{\|\delta f\|_2^2}}
\ee
for all $f$ and $t>0$. On the other hand, we have the Poincar\'{e} inequality: there exists
some $c_2>0$ such that
\be\label{poincare}
\E\left[ (f-\E f)^2\right]\leq c_2 \sum_{i\in\Z} \int (\nabla_i f)^2 \dd\pee
\ee
for all $f$.
Now choosing
\[
f = \caA_N^\ell
\]
we easily see that
\[
\|\delta f\|_2^2 \leq \ell^2/N\;.
\]
This combined with \eqref{gauss} gives \eqref{colibri}.
To see that this implies almost-sure convergence to $\E(\si_0)^\ell$, we use the strong
mixing property enjoyed by one-dimensional Gibbs measures with exponentially decaying interacting
\cite[Chap. 8]{geo}, from which it follows easily that
\[
|\E(\si_{ki}|\si_{ri}, r\not= k)-\E(\si_0)|\leq Ce^{-ci}
\]
which implies
\[
|\E (\si_{i}\si_{2i}\cdots \si_{\ell i})-\E(\si_0)^\ell |\leq C_\ell\, e^{-ci}
\]
This in turn implies 
\[
\lim_{N\to\infty} \E(\caA^\ell_N)= \E(\si_0)^\ell.
\]
Combining this fact with \eqref{colibri} yields the almost-sure convergence
of   $\caA^\ell_n$ towards $\E(\si_0)^\ell$ as $n$ goes to infinity.
The first statement is thus proved.

In order to prove the second statement,
we use the bound
\be\label{bim}
\E\left(\prod_{j=1}^q \si_{i_j}\right)\leq e^{-\gamma q}
\ee
for some $\gamma>0$ and for all $i_1,\ldots,i_q\in \Z$.
This follows immediately from the `finite-energy property' of one-dimensional
Gibbs measures, i.e., the fact that there exists $\delta\in (0,1)$ such that
for all $\si\in\Omega, \alpha\in \{0,1\}$
\[
\delta<\pee \left(\si_0=\alpha|\si_{\Z\setminus\{0\}}\right)<1-\delta.
\] 
As a consequence,
\[
\left|\nabla_j \left(\prod_{r=1}^\ell \si_{ir}\right)\right|\leq \1(j\in \{i,2i,\ldots,\ell i\})
\prod_{r=1, ri\not=j}^\ell \si_{ri}
\]
and hence, using the elementary inequality $(\sum_{i=1}^N a_i)^2\leq N\sum_{i=1}^n a_i^2$,
we have the upper bound
\[
\left|\nabla_j \caK(N,\ell)\right|^2
\leq
N\sum_{i=1}^{\lfloor N/\ell\rfloor} \prod_{r=1}^\ell \si_{ir} \1(j\in \{i,2i,\ldots,\ell i\}).
\]
Integrating against $\pee$, using \eqref{bim} and summing over $j$ yields
\[
\sum_{j} \int \left(\nabla_j \caK(N,\ell)\right)^2 \dd\pee\leq N^2 e^{-\ell\gamma}\;.
\]
Choosing now
\[
\ell=\ell(N)= c\log N,
\]
and using \eqref{poincare}, we find
\[
\textup{Var} (\caK(N,\ell(N)))\leq CN^{2-c\gamma}\;.
\]
Hence, for $\gamma>2/c$, the variance of $\caK(N,c\log N)$ converges to zero. Since
\[
\E(\caK(N,\ell(N))\leq N\, e^{-\gamma k}\leq N^{1-c\gamma},
\]
the expectation of $\caK(N,\ell(N))$  also converges to zero, hence we have convergence to zero in mean
square sense and thus in probability.
\epr
\section{Large deviations for arithmetic progressions of size two}\label{sect-size2}

From the point of view of functional inequalities such as the Gaussian
concentration bound or the Poincar\'{e} inequality, there is hardly a
difference between sums of shifts of a local function, i.e.
conventional ergodic averages, and their nonconventional counterparts.

The difference becomes however manifest in the study of large deviations.
If we think e.g.\ about $\sum_{i=1}^N\si_{i}\si_{i+1}$
versus $\sum_{i=1}^N \si_i \si_{2i}$
as ``Hamiltonians'' then the first sum
is simply a nearest neighbor translation-invariant interaction,
whereas the second sum is a long-range non translation
invariant interaction.
Therefore, from the point of view of computing partition
functions, the second Hamiltonian will be much
harder to deal with.

In this section we restrict to the product case, by
choosing $\pee_p$ to be product of Bernoulli
with parameter $p$ on two symbols $\{ +,-\}$, and
consider arithmetic progressions of size two $(k=2)$.
We will show that the thermodynamic limit of the
free energy function associated to the sum
\[
S_N= \sum_{i=1}^N \si_i\si_{2i}
\]
defined
as
\be\label{bom}
\caF_p(\lambda)=\lim_{N\to\infty}\frac1N\log\E_p\left(e^{\lambda S_N}\right)
\ee
exists, is analytic as a function of $\lambda$ and
has an explicit expression in terms of combinations of Ising model partition functions
for different volumes.

To start, assuming $N$ to be odd (the case $N$ even is treated similarly),
we make the following useful decomposition
\[
S_N= \sum_{l=1}^{\frac{N+1}{2}} S_l^{(N)}
\]
with
\be\label{kok}
S_l^{(N)}= \sum_{i=0}^{M_l(N)-1} \si_{(2l-1)2^i}\si_{(2l-1)2^{i+1}}
\ee
and
\[
M_l(N)=\left\lfloor \log_2 \left(\frac{N}{2l-1}\right)\right\rfloor +1\;.
\]
where $\lfloor x\rfloor$ denotes the integer part of $x$.
The utility of such decomposition is that the random variable
$S_l^{(N)}$ is independent
from $S_{l'}^{(N)}$ for $l\neq l'$. A  similar decomposition into independent blocks has also been used
independently in \cite{FLM,FSW}.
This implies that the partition function
in the free energy \eqref{bom} will factorize over different subsystems
labeled by $l\in\{1,\ldots, (N+1)/2\}$, each of size $M_l(N) +1$.
Therefore we can treat separately each variable $S_l^{(N)}$.

Furthermore, defining new spins
\[
\tau_i^{(l)} = \sigma_{(2l-1)2^{i-1}} \qquad \text{for} \quad i \in \{1,\ldots, M_l(N)+1\}\;,
\]
it is easy to realize that, for a given $l\in\{1,\ldots, (N+1)/2\}$,  the variable
$S_l^{(N)}$ is nothing else than the Hamiltonian
of a one-dimensional nearest-neighbors Ising model, since
\[
\{S_l^{(N)}\} \stackrel{\EuScript{D}}{=} \left\{\sum_{i=1}^{M_l(N)}\tau^{(l)}_i\tau^{(l)}_{i+1}\right\}
\]
where $\tau^{(l)}_i$ are Bernoulli random variables with parameter $p$,
independent for different values of $l$ and for different values of $i$ and $\stackrel{\EuScript{D}}{=}$
denotes equality in distribution.
Introduce the notation
\[
\caZ(\lambda,h,n+1)= \sum_{\tau\in \{-1,1\}^{n+1}} e^{\lambda\sum_{i=1}^n \tau_i\tau_{i+1} + h\sum_{i=1}^{n+1}\tau_i}
\]
for
the partition function of the one-dimensional
Ising model with coupling strength $\lambda$ and external field $h$ in the volume
$\{1,\ldots, n\}$, with {\em free} boundary conditions.
Then we have
\be\label{mercredi}
\E_p \left(e^{\lambda\sum_{i=1}^n \tau_i\tau_{i+1}}\right)= (p(1-p))^{\frac{n+1}{2}}\caZ (\lambda,h,n+1)
\ee
with $h= \frac12\log(p/(1-p))$. A standard computation 
(see for instance \cite{baxter1982exactly}, Chapter 2) 
gives
\[
\caZ(\lambda,h,n+1)= v^T M^n v= | v^T \cdot e_+|^2 \Lambda_+^n
+ |v^T\cdot e_-|^2 \Lambda_-^n
\]
with
$\Lambda_{\pm}$ the largest, resp.\ smallest eigenvalue of the transfer matrix
(with elements $M_{\alpha, \beta}= e^{\lambda \alpha\beta + \frac{h}{2}(\alpha+\beta)}$), i.e.,
\[
\Lambda_{\pm}= e^\lambda \left(\cosh(h)\pm \sqrt{\sinh^2(h) + e^{-4\lambda}}\right),
\]
$v^T$ the vector with components $(e^{h/2}, e^{-h/2})$, $e_{\pm}$ the normalized eigenvectors corresponding to the eigenvalues
$\Lambda_{\pm}$.

Using the decomposition \eqref{kok}, we obtain from \eqref{mercredi}
\[
\log\E_p \left(e^{\lambda S_N}\right)=
\sum_{l=1}^{(N+1)/2}\log\left(p(1-p)^{\frac{M_l(N)+1}2}\caZ (\lambda,h,M_l(N)+1)\right).
\]
Furthermore, observing that
\[
\lim_{N\to\infty}\frac1N \sum_{l=1}^{(N+1)/2} M_l(N)= \frac 1 2 \int \psi(x) \dx
\]
with
\[
\psi(x)= \left\lfloor \log_2\left(\frac1x\right)\right\rfloor +1\;,
\]
we obtain
\begin{align*}
\caF_p (\lambda) &= \frac1{4}\left(\int \psi(x) \dx +1\right)  \log (p(1-p))\nonumber\\
& \quad + \frac12 \int_0^1 \log \left( |v^T \cdot e_+|^2 \Lambda_+^{\psi(x)} +
|v^T \cdot e_-|^2 \Lambda_-^{\psi(x)}\right)\dx.
\end{align*}
To obtain a more explicit formula one can make use of the following:
the normalized eigenvector corresponding to the largest eigenvalue is
\[
e_+ =\frac{w_+}{\|w_+\|}
\]
with
\[
w_+ =
\left(
\begin{array}{c}
-e^{-\lambda} \\
e^{h+\lambda}-\Lambda_+
\end{array}
\right)
\]
and moreover
\[
|v^T \cdot e_-|^2= \|v\|^2- |v^T \cdot e_+|^2 = 2\cosh (h) - |v^T \cdot e_+|^2 \;.
\]
Since $\psi(x)= n+1$ for $x\in (1/2^{n+1}, 1/2^n]$, we have
\[
\int \psi(x) \dx = 2 \;,
\]
hence one gets
\be\label{fexpli}
\caF_p(\lambda) = \log\left([p(1-p)]^\frac 3 4 |v^T \cdot e_+|\; \Lambda_+\right) + \caG(\lambda)
\ee
with
\[
\caG(\lambda)= \frac12 \sum_{n=1}^\infty \frac{1}{2^{n}} \log\left(1+ \left(\frac{2\cosh(h)}{|v^T \cdot e_+|^2}-1\right)\left(\frac{\Lambda_-}{\Lambda_+}\right)^{n}\right)\;.
\]
In the case $p=1/2$, we have $h=0$, $\Lambda_+= e^\lambda+ e^{-\lambda}$, $|v^T \cdot e_+|^2= \|v\|^2=2$ which implies
$\caG(\lambda)=0$ and
\be\label{iban}
\caF_{1/2}(\lambda) = \log\left(\frac12\left(e^{\lambda}+ e^{-\lambda}\right)\right)\;.
\ee
One recognizes in this case the Legendre transform of the large deviation rate function for a sum of i.i.d. bernoulli $(1/2)$ because (only) in this case $p=1/2$
the joint distribution of $\{\si_i \si_{2i}, i\in \N\}$ coincides with the joint distribution
of a sequence of independent Bernoulli$(1/2)$ variables.
When $p\not=1/2$, although an explicit formula is
given in \eqref{fexpli}, the expression reflects the multiscale
character of the decomposition and it is non-trivial.

As a consequence of the explicit formula
\eqref{fexpli},
we have the following

\bt\label{th1}\par\noindent
\ben
\item {\em Large deviations.} The sequence of random variables
$\frac{S_N}{N}$ satisfies a large deviation principle
with rate function
\[
I_p(x)= \sup_{\lambda\in\R} \left(\lambda x- \caF_p(\lambda)\right)
\]
where $\caF_p(\lambda)$ is given by \eqref{fexpli}.
\item {\em Central limit theorem.} The sequence of random variables
\[
N^{-1/2}\left(S_N-\E_p(S_N)\right)
\]
weakly converges to a Gaussian random variable with strictly positive
variance $\si^2= \caF_p''(0)>0$. 
\een
\et
\bpr
The expression \eqref{fexpli}
shows that $\caF_p$ is differentiable as a function of $\lambda$, hence
the first statement follows from the G\"{a}rtner-Ellis theorem \cite{dz}.
The second statement follows from the fact that $\caF_p$ is analytic
in a neighborhood of the origin, which again follows
directly from the explicit expression.
\epr

\br
The value $p=1/2$ is special since in this case the joint distribution of $\{\si_i \si_{2i}, i\in \N\}$
coincides with the joint distribution of a sequence of independent $\textup{Bernoulli}(1/2)$ variables.
Therefore we must have
\[
I_{1/2}(x)= 
\begin{cases}
\frac12 (1+x)\log\left(1+x\right) +
\frac12 (1-x)\log\left(1-x\right) & \text{ if } |x|\le 1\;,\\
+\infty & \text{ if } |x|>1\;.
\end{cases} 
\]
This can be checked by computing $\caF_{1/2}$ as done above (see \eqref{iban}).
We must also have $\si^2=1$ for $p=1/2$.
\er

\br
Notice that we computed the large deviation rate function in the  $\pm 1$ setting.
If one considers a Bernoulli measure ${\mathbb Q}_p$ on $\{0,1\}^{\Z}$,
with ${\mathbb Q}_p (\eta_i=1)=p$, then the large deviations of the sums
\be\label{boemboem}
\sum_{i=1}^N \eta_i\eta_{2i}
\ee
correspond to the large deviations of 
\[
\frac14 \sum_{i=1}^N (1+\si_i)(1+\si_{2i})=
\frac14 \left(N + \sum_{i=1}^N (\si_i+\si_{2i}) + \sum_{i=1}^N \si_i\si_{2i}\right)
\]
where $\si$ is distributed according to $\pee_p$ on $\{+,-\}^{\Z}$. In particular, the free energy for the 
large deviations of \eqref{boemboem} under the measure ${\mathbb Q}_{1/2}$ corresponds to a free energy
of the $\si$ spins with non-zero magnetic field and hence can again be computed explicitly.
\er

\begin{figure}[htbp!!!!]
\centering
\includegraphics[scale=.7]{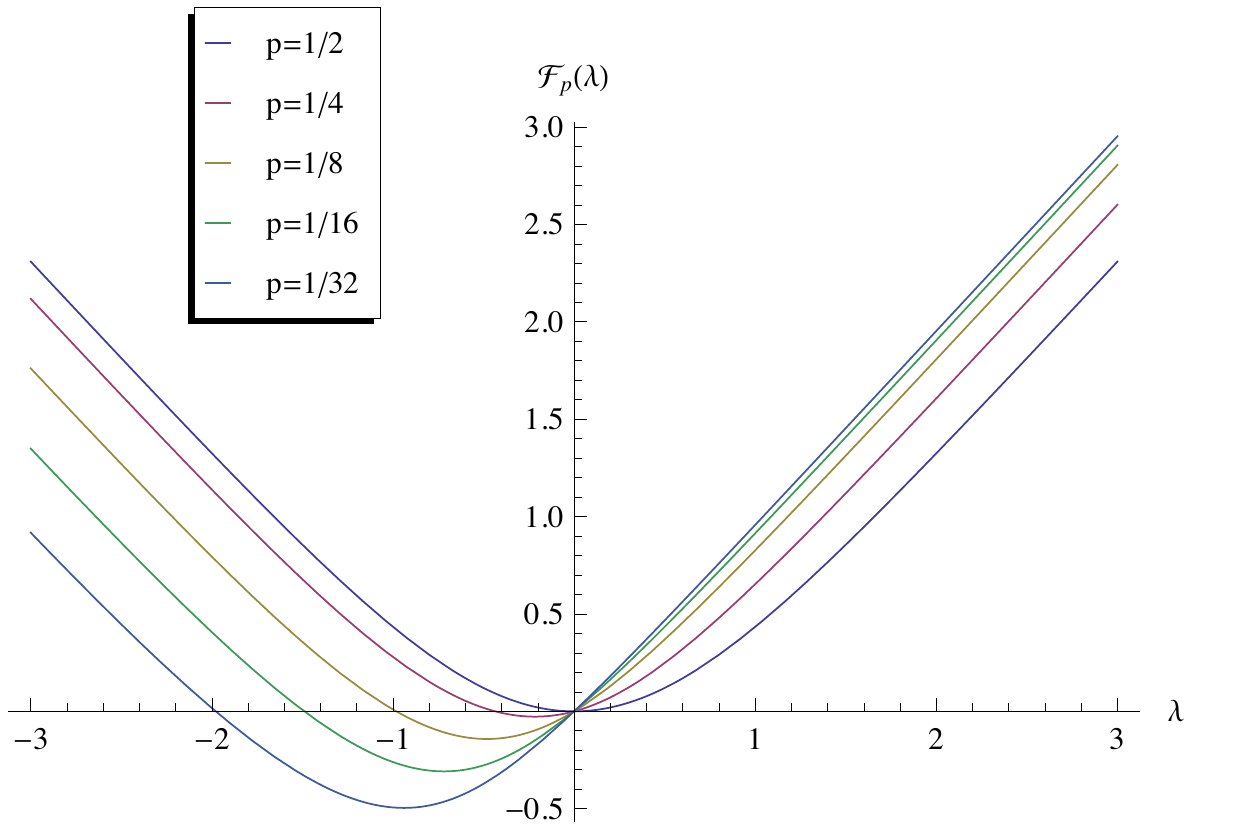}
\caption{
{\small Plot of the free energy function for different $p$ values.
The graph has been obtained from formula (\ref{fexpli})
truncating the sum to the first 100 terms.}}
\label{FreeEng}
\end{figure}

\begin{figure}[h!!!!]
\centering
\includegraphics[scale=.7]{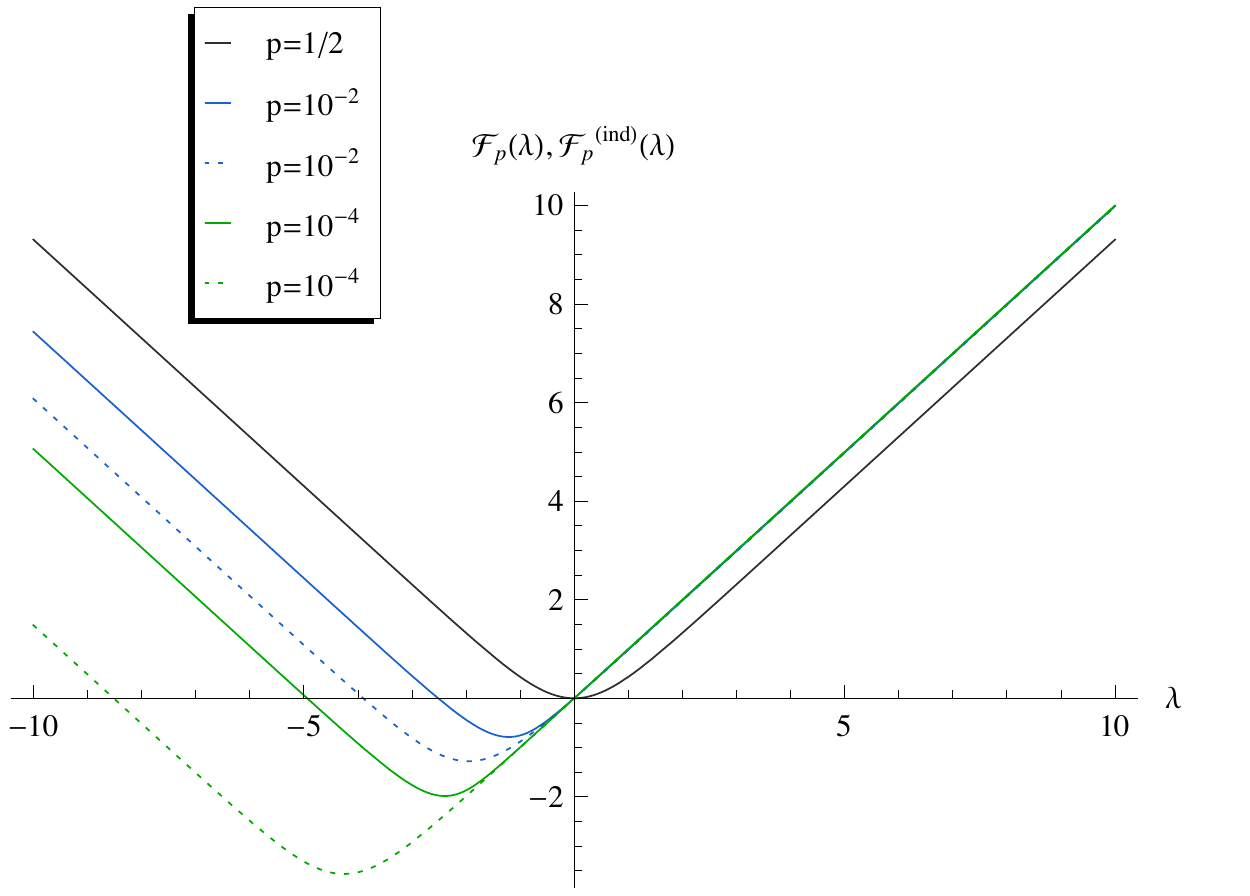}
\caption{
{\small Plot of the free energy function (\ref{fexpli}) (continuous line)
and of the free energy of the independent case (\ref{inde}) (dashed line).
}}
\label{FreeEng5}
\end{figure}

\br
A plot of the free energy for a few values of $p$ is shown in
Figure \ref{FreeEng} (it is enough to analyze values in $(0,1/2]$
since $\caF_{p}(\lambda) = \caF_{1-p}(\lambda)$).
In the general case $p\neq 1/2$ it is interesting to compare our results
to the independent case. To this aim one consider the sum
$\sum_{i=1}^N \xi_i\eta_i$ where $\xi_i, \eta_i$ are two sequences
of i.i.d. Bernoulli of parameter $p$. Note that in this case
the family $\{\xi_i \eta_i\}_{i\in\{1,\ldots,N\}}$ is made of
independent Bernoulli random variables with parameter
$p^2 + (1-p)^2$. An immediate computation of the
free energy yields on this case
\be
\label{inde}
 \mathcal F_p^{(\textup{ind})}(\lambda)= \log \left([p^2 + (1-p)^2]e^{\lambda} + 2p(1-p) e^{-\lambda}\right)
\ee
This free energy is compared to that of formula (\ref{fexpli}) in
Figure \ref{FreeEng5}.
\er

\begin{figure}[htbp!!!]
\centering
\includegraphics[scale=.7]{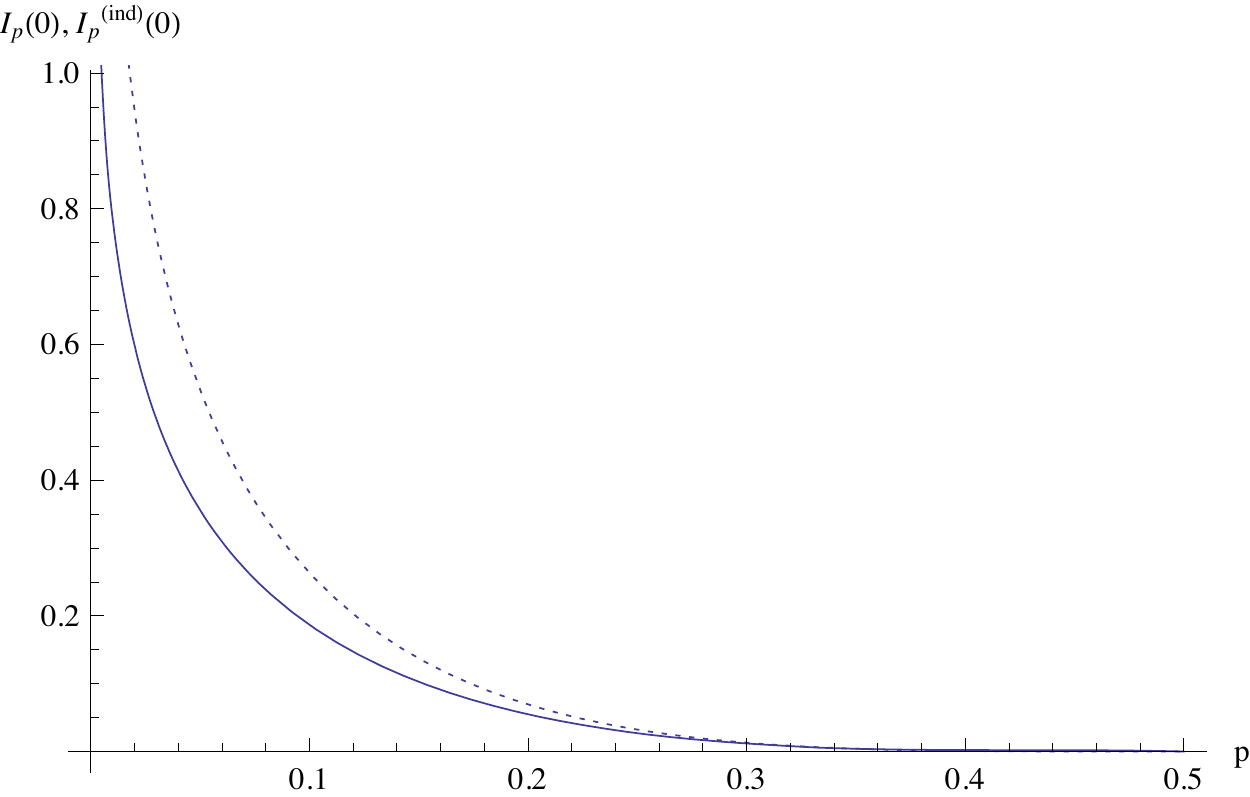}
\caption{
{\small The negative of the minimum of the free energy as a
function of $p$ for the two cases (\ref{fexpli}) and (\ref{inde}).
}}\label{FreeEngMin2}
\end{figure}

In particular one can analyze the behaviour of the minimum
of the free energy functions in the two cases, corresponding
to the negative value of the large deviation rate function computed
at zero. This is shown in Figure \ref{FreeEngMin2}, which suggests
a general inequality between the two cases.

\newpage

\section{Size larger than  two and Ising model in higher dimension}\label{cloclo}

In this Section we analyze the case of arithmetic progressions
of size larger than two. Such a case naturally leads to statistical 
mechanics models in dimension higher than one, possibly having
phenomenon of phase transitions. Conversely the classical
Ising model in dimension $d>1$ can be related to specific
unconventional sums, which we describe below.

\subsection{Decompositions for $k \ge 2$}

When the size of the arithmetic progressions is larger than two ($k>2$),
we have sums of the type
\[
S_N =\sum_{i=1}^N \sigma_i \sigma_{2i} \cdots \sigma_{ki}\;,
\]
where $\sigma_i$ are i.i.d. random variables taking values
in the set $\{-1,+1\}$.
One can try to decompose this sum into independent sums
as it was done in Section \ref{sect-size2}.
After relabeling  the indices one obtains independent sums,
each of which corresponds to a spin system with Hamiltonian
in a bounded domain of $\N^{d_k}$,
where the dimension $d_k$ 
is given by the number of prime numbers contained in the set $\{2, \ldots, k\}$.
Denoting by  $p_1, p_2,\ldots, p_{d_k}$   the prime numbers contained in $\{2, \ldots, k\}$
and defining
\[
A^{(N)}_{p}:= \{m \in\{2,\ldots,N\}\,: \; m \: \text{ is not divisible by }\:  p\}
\]
and
\[ 
A^{(N)}_{p_1,\ldots,p_{d_k}}:= \bigcap_{l=1}^{d_k} A_{p_l}^{(N)}\;
\]
then one has the following decomposition:
\[
S_N = \sum_{m\in A^N_{p_1,\ldots,p_{d_k}} } S_m^{(N)}\;.
\]
The independent sums $S_m^{(N)}$ are   given by 
translation-invariant Hamiltonians of the form
\[
S_m^{(N)} =  \sum_{X\subset \N^{d_k}} J_X^{(m)} \tau_X^{(m)}
\]
where the spins $\tau_X^{(m)}$ are given by 
\[
\tau_X^{(m)} = \prod_{{j}\in X} \tau_{{j}}^{(m)}
:= \prod_{(j_1, j_2, \ldots, j_{d_k})\in X} \sigma_{m \, p_1^{j_1}\, p_2^{j_2}\cdots p_{d_k}^{j_k} }
\]
and the couplings $J_X^{(m)}$ are
\[
J_X^{(m)} = 
\begin{cases}
1 & \text{if } X = T_{l} X(k) \quad \text{for some} \quad {l}\in \Lambda_{p_1,\ldots,p_{d_k}}^m(N)\\
0 & \text{otherwise}
\end{cases} 
\]
with
\[
\Lambda_{p_1,\ldots,p_{d_k}}^m(N)=\left\{0,1,\ldots,\left\lfloor
\log_{p_1} \frac N m \right\rfloor \right\}\times \cdots\times \left\{0, 1,\ldots, \left\lfloor\log_{p_{d_k}}
\frac N m\right\rfloor \right\} \subset \mathbb{N}^{d_k}\,,
\]
$T_l X$ the translation of the set $X$ by the vector $l$, 
and $X(k)$ a specific subset of $\N^{d_k}$ depending on the
size of the arithmetic progression $k$.
This set $X(k)$ is a polymer starting at the origin and  having  $k$ vertices.
The specific shape of $X(k)$ sets the range of interaction along 
each direction of the $d_k-$dimensional lattice. 
In general the shape of the interaction  depends
on the non-prime numbers contained in $\{2,\ldots,k\}$.\\

We clarify this construction with a few examples.

\begin{itemize}

\item $k=2$
\[
\sum_{i=1}^N \sigma_i \sigma_{2i}= \sum_{m \in A_2^{(N)}} \sum_{i\in \Lambda_2^m(N)} \sigma_{m \cdot 2^i} \:\sigma_{m \cdot 2^{i+1}}= \sum_{m \in A_2^{(N)}} \sum_{i\in \Lambda_2^m(N)} \tau_i^{(m)}  \tau_{i+1}^{(m)}\;.
\]
This Hamiltonian is the 1-dimensional nearest-neighbor Ising model
of Section \ref{sect-size2} constructed from the basic polymer $X(2) = \{0,1\}$.

\item $k=3$
\begin{align*}
\sum_{i=1}^N \sigma_i \sigma_{2i} \sigma_{3i}
&= \sum_{m \in A_{2,3}^{(N)} } \sum_{(i,j)\in \Lambda_{2,3}^m(N)} \sigma_{m \cdot 2^{i} 3^{j}} \:\sigma_{m \cdot 2^{i+1} 3^{j}}\: \sigma_{m \cdot 2^{i} 3^{j+1}} \nonumber\\
&= \sum_{m \in A_{2,3}^{(N)} } \sum_{(i,j)\in \Lambda_{2,3}^m(N)}  \tau_{i,j}^{(m)} \, \tau_{i+1,j}^{(m)} \,\tau_{i,j+1}^{(m)}  \;.
\end{align*}

This corresponds to a 2-dimensional nearest-neighbor model
with triple interaction obtained via the polymer $X(3) = \{(0,0),(0,1),(1,0)\}$.

\item $k=4$
\begin{align*}
\sum_{i=1}^N \sigma_i \sigma_{2i} \sigma_{3i} \sigma_{4i}
&= \sum_{m \in A_{2,3}^{(N)} } \sum_{(i,j)\in \Lambda_{2,3}^m(N)} \:\sigma_{m \cdot 2^{i} 3^{j}} \:\sigma_{m \cdot 2^{i+1} 3^{j}} \:\sigma_{m \cdot 2^{i+2} 3^{j}}\: \sigma_{m \cdot 2^{i} 3^{j+1}} \nonumber\\
&= \sum_{m \in A_{2,3}^{(N)} } \sum_{(i,j)\in \Lambda_{2,3}^m(N)}\tau_{i,j}^{(m)}  \,\tau_{i+1,j}^{(m)}\, \tau_{i+2,j}^{(m)} \,\tau_{i,j+1}^{(m)} \;. 
\end{align*}

This gives a 2-dimensional model sums with quadruple interaction constructed 
by translating the polymer $X(4) = \{(0,0),(1,0),(2,0),(0,1)\}$.
The range of interaction is 2 in one direction and 1 in the other direction.

\item $k=5$
\begin{align*}
&\sum_{i=1}^N \sigma_i \sigma_{2i} \sigma_{3i} \sigma_{4i} \sigma_{5i}\nonumber\\
&=\sum_{m \in A_{2,3,5}^{(N)} } \sum_{(i,j,l)\in \Lambda_{2,3,5}^m(N)}\sigma_{m \cdot 2^i 3^j 5^l} 
\:\sigma_{m \cdot 2^{i+1} 3^j 5^l} \:\sigma_{m \cdot 2^{i+2} 3^j 5^l}\: \sigma_{m \cdot 2^i 3^{j+1} 5^l} \: \sigma_{m \cdot 2^i 3^{j} 5^{l+1}} \nonumber\\
&= \sum_{m \in A_{2,3,5}^{(N)} } \sum_{(i,j,l)\in \Lambda_{2,3,5}^m(N)} \tau_{i,j,l}^{(m)} \, \tau_{i+1,j,l}^{(m)} \,\tau_{i+2,j,l}^{(m)} \,\tau_{i,j+1,l}^{(m)} \, \tau_{i,j,l+1}^{(m)}\;.
\end{align*}
Here we get a 3-dimensional model with quintuple interaction
given by the basic polymer $X(5) = \{(0,0,0),(1,0,0),(2,0,0),(0,1,0),(0,0,1)\}$. The range of interaction is 2 in one direction and 1 in the other two directions.

\end{itemize}

\subsection{Unconventional sums related to 2-dimensional Ising model}

We consider now the standard 2-dimensional neirest-neighbor Ising model sums
\[
 \sum_{(i,j)\in\La}\tau_{i,j} ( \tau_{i+1,j} +\tau_{i,j+1})
\]
in a domain $\La$ of $\N^2$, and wonder whether there exist  some unconventional averages that may be related to it through the decomposition procedure previously described. The answer is in the affirmative sense and is contained in the following
two examples.

\begin{itemize}

\item
For $\{\sigma_i\}_{i \in \N}$ a sequence of independent random variables
taking values in $\{-1,+1\}$, we have 
\begin{align*}
\sum_{i=1}^N \sigma_i (\sigma_{2i} +\sigma_{3i})
&= \sum_{m \in A_{2,3}^{(N)}}  \sum_{(i,j)\in \Lambda_{2,3}^m(N)}   \sigma_{m \cdot 2^i 3^j} \:(\sigma_{m \cdot 2^{i+1} 3^j}+ \sigma_{m \cdot 2^i 3^{j+1}} )\nonumber\\
&= \sum_{m \in A_{2,3}^{(N)}}  \sum_{(i,j)\in \Lambda_{2,3}^m(N)} \tau_{i,j}^{(m)} \, (\tau_{i+1,j}^{(m)} +\tau_{i,j+1}^{(m)})
\end{align*}
with $\tau_{i,j}^{(m)}= \sigma_{m \cdot 2^i 3^j}$.
This clearly gives a decomposition into $|A_{2,3}^{(N)}|$ independent two-dimensional nearest-neighbor Ising sums.
\item
Let $\sigma_{i,j}$ be i.i.d. dichotomic random variables labeled by $(i,j)\in \N^2$. Then
\begin{align*}
\sum_{i,j=1}^N \sigma_{i,j} (\sigma_{2i,j} +\sigma_{i,2j})
&= \sum_{m \in A_{2}^{(N)}}  \sum_{(i,j)\in \Lambda_{2}^m(N)}   \sigma_{m \cdot 2^i, m \cdot 2^j} \:\left(\sigma_{m \cdot 2^{i+1}, m \cdot 2^j}+ \sigma_{m \cdot 2^i, m \cdot 2^{j+1}} \right)\nonumber\\
&= \sum_{m \in A_{2}^{(N)}}  \sum_{(i,j)\in \Lambda_{2}^m(N)} \nu_{i,j}^{(m)} \, (\nu_{i+1,j}^{(m)} +\nu_{i,j+1}^{(m)})
\end{align*}
with $\nu_{i,j}^{(m)}:=\sigma_{m \cdot 2^i, m \cdot 2^j}$. We have a decomposition into $|A_{2}^{(N)}|$ independent two-dimensional nearest-neighbor Ising sums.

\end{itemize}


\end{document}